\documentclass[stef,twoside]{stefano}

\usepackage{amsmath}
\usepackage{amssymb}
\usepackage{graphics}
\usepackage{rotating}
\usepackage{cite}
\usepackage{color}

\def\makeheadbox{{%
\hbox to0pt{\vbox{\baselineskip=10dd\hrule\hbox
to\hsize{\vrule\kern3pt\vbox{\kern3pt \hbox{  {\sc El. J. Lin. Alg.} {\bf
15}, 297-313 (2006)} \hbox{ {\sc
{\color{blue}{dma}}[{\color{black}{imecc}}]{\color{red}{UniCamp}} }
\hspace*{10.4cm} {\color{blue}{$\boldsymbol{\Sigma \delta \Lambda}$}} }
\kern3pt}\hfil\kern3pt\vrule}\hrule}%
\hss}}}

\def\0{\mbox{\tiny $0$}}
\def\1{\mbox{\tiny $1$}}
\def\2{\mbox{\tiny $2$}}
\def\3{\mbox{\tiny $3$}}
\def\4{\mbox{\tiny $4$}}
\def\5{\mbox{\tiny $5$}}
\def\6{\mbox{\tiny $6$}}
\def\7{\mbox{\tiny $7$}}
\def\8{\mbox{\tiny $8$}}
\def\9{\mbox{\tiny $9$}}

\def\mi{\mbox{\tiny $-$}}

\def\eq{\mbox{\tiny $=$}}

\begin{document}

\title{\Large ZEROS OF UNILATERAL QUATERNIONIC POLYNOMIALS}

\author{
Stefano De Leo\inst{1},
 Gisele Ducati\inst{2},
\and Vinicius Leonardi\inst{3} }

\institute{
Department of Applied Mathematics, University of Campinas\\
PO Box 6065, SP 13083-970, Campinas, Brazil\\
{\em deleo@ime.unicamp.br}\and
Department of Mathematics, University of Parana\\
PO Box 19081, PR 81531-970, Curitiba, Brazil\\
{\em ducati@mat.ufpr.br}  \and
Department of Physics, University of Parana\\
PR 81531-970, Curitiba, Brazil\\
{\em vjhcl02@fisica.ufpr.br}
}


\date{Submitted: {\em November, 2005} - Revised: {\em Septenber, 2006}.}










\abstract{The purpose of this paper is to show how the problem of finding
the zeros of  unilateral $n$-order quaternionic polynomials can be solved
by determining the eigenvectors of the corresponding companion matrix. This
approach, probably superfluous in the case of quadratic equations for which
a closed formula can be given, becomes truly useful for (unilateral)
$n$-order polynomials. To understand the strehgth of this method, we
compare it with the Niven algorithm and show where this (full) matrix
approach improves previous methods based on the use of the Niven algorithm.
For the convenience of the readers, we explicitly solve some examples of
second and third order unilateral quaternionic polynomials. The leading
idea of the practical solution method proposed in this work can be
summarized  in following three steps: translating the quaternionic
polynomial in
 the eigenvalue problem for its companion matrix, finding its eigenvectors,
 and, finally, giving the quaternionic solution of the unilateral polynomial
 in terms of the components of such eigenvectors. A brief
 discussion on bilateral quaternionic quadratic equations is also
 presented.}

\authorrunning{S. De Leo, G. Ducati, V. Leonardi}
\titlerunning{\sc unilateral quaternionic polynomials}

\maketitle


\section*{I. INTRODUCTION}

As far as we know, the problem of solving quaternionic quadratic equations
and the study of the fundamental theorem of algebra for quaternions were
first approached by Niven and Eilenberg \cite{NIV42,EN44}. After these
fundamental works, the question of counting  the number of zeros of
$n$-order quaternionic polynomials and the problem to find the solutions
have been investigated. An interesting application of unilateral quadratic
equations (where a closed formula for their solutions can be obtained) is
found in solving homogeneous quaternionic linear second order differential
equations with constant coefficients \cite{DD01,DD03}. The solution of
\begin{equation}
\label{2ode} \frac{\mbox{d}^{\2}\,\,}{\mbox{d}x^{\2}} \, \Psi(x) -
a_{\1}\,\frac{\mbox{d}\,\,}{\mbox{d}x} \,\Psi(x) - a_{\0}\,\Psi(x)
= 0\,\,,
\end{equation}
where $\Psi\,:\,\mathbb{R}\to \mathbb{H}\,\,,\,\,a_{\0,\1} \in \mathbb{H}$
(the set of quaternions) and $x\in \mathbb{R}$, can indeed be reduced, by
setting $\Psi(x)=\exp[q\,x]$ ($q\in\mathbb{H}$) and using the
$\mathbb{H}$-linearity of Eq.(\ref{2ode}), to the following quadratic
equation
\begin{equation}
\label{p2} q^{\2} - a_{\1}\,q - a_{\0}=0\,\,.
\end{equation}
Generalizing this result, we can see that the problem of finding the
solution of $n$-order ($\mathbb{H}$-linear) differential equations with
constant coefficients,
\begin{equation}
\label{node} \frac{\mbox{d}^n\,\,}{\mbox{d}x^n} \,\Psi(x) -\sum_{s
\eq \1}^{n \mi\1} a_{s}\, \frac{\mbox{d}^s\,\,}{\mbox{d}x^s}
\,\Psi(x) - a_{\0}\,\Psi(x)=0\,\,,
\end{equation}
where $a_{\0,\1,...,n \mi \1} \in \mathbb{H}$, can be immediately
transformed to the problem of finding the zeros of the corresponding
$n$-order unilateral quaternionic polynomial,
\begin{equation}
\label{pn} A_n(q) :=  q^n - \sum_{s \eq \0}^{n \mi\1} a_{s}\,
q^{s}\,\,.
\end{equation}
The important point to be noted here is that such a solution
method based on  finding the zeros of the corresponding
quaternionic polynomial for  $\mathbb{H}$-linear differential
equations with constant coefficients  does not work for
$\mathbb{R}$ and $\mathbb{C}$-linear differential equations. For
example, in non-relativistic quaternionic quantum
mechanics\cite{ADL}, the one-dimensional Schr\"odinger equation
(for quaternionic stationary states) in presence of a constant
quaternionic potential, reads
\begin{equation}
\label{se} i\,\frac{\hbar^{\2}}{2m}\,
\frac{\mbox{d}^{\2}\,\,}{\mbox{d}x^{\2}} \, \, \Phi(x)-
\left(i\,V_{\1}+j\,V_{\2}+k\,V_{\3}\right)\,\Phi(x) = - \,
\Phi(x)\,E\, i\,\,,
\end{equation}
and, due to the $\mathbb{C}$-linearity, its solution cannot be
expressed in terms of a quaternionic exponential function,
$\exp[q\,x]$. Recent studies of quaternionic barrier \cite{DDN02}
and well \cite{DD05} potentials confirmed that the solution of
Eq.(\ref{se}) has to be expressed as the product of two factors: a
quaternionic constant coefficient, $p$, and a complex exponential
function, $\exp[z\,x]$ ($z\in \mathbb{C}$). Using
$\Phi(x)=p\,\exp[z\,x]$, the Eq.(\ref{se}) reduces to
\begin{equation}
\label{red} i\,\frac{\hbar^{\2}}{2m}\, p\,z^{\2} -
\left(i\,V_{\1}+j\,V_{\2}+k\,V_{\3}\right)\,p  = - \, p \,E\,
i\,\,,
\end{equation}
which obviously does not represent a quaternionic unilateral
polynomial. For a more complete discussion of $\mathbb{R}$ and
$\mathbb{C}$ linear quaternionic differential equation theory, we
refer the reader to the papers cited in refs.\cite{DD01,DD04}.

The problem of finding the zeros of unilateral quaternionic polynomials was
solved in 1941 by using a  {\em two-step} algorithm. In his seminal work
\cite{NIV41}, Niven proposed to divide the unilateral $n$-order quaternion
polynomial (\ref{pn})  by a quadratic polynomial with real coefficients
($c_{\0,\1}\in \mathbb{R}$)
\begin{equation}
\label{c2} C_{\2}(q) := q^{\2} - c_{\1}\,q - c_{\0}\,\,.
\end{equation}
Then, by using the polynomial equation
\begin{equation}
\label{pe} A_n(q) = B_{n \mi \2}(q) \, C_{\2}(q) - D_{\1}(q)\,\,,
\end{equation}
where
\[ B_{n \mi \2}(q) :=  q^{n\mi \2} - \sum_{s
\eq \0}^{n \mi\3} b_{s}\, q^{s}\,\,,\,\,\,\,\, D_{\1}(q) :=
d_{\1}\,q + d_{\0}\,\,,\,\,\, \,\,b_{\0,\1,...n \mi \3}\,,\,\,
d_{\0,\1} \in \mathbb{H}\,\,,
\]
 the problem of finding the zeros of $A_n(q)$ was translated
in the following {\em two-step} problem:
\begin{itemize}
\item step 1  -  determine  $d_{\0}$ and $d_{\1}$ in terms of $c_{\0,\1}$
and $a_{\0,\1,...,n \mi \1}$;
\item step 2 - obtain two coupled real
equations which allow to calculate $c_{\0}$ and $c_{\1}$.
\end{itemize}
To facilitate the understanding of this paper, we explicitly discuss the
Niven algorithm in the next section. Since its proof is quite detailed and,
furthermore, since its use can be {\em completely} avoided, we merely state
the results without proof.

In a recent paper \cite{SER01}, Ser\^odio, Pereira and Vitoria (SPV)
proposed  a practical technique to compute the zeros of unilateral
quaternion polynomials. The SPV approach is based on the use of the Niven
algorithm \cite{NIV41}. Their idea was to improve such an algorithm by
calculating the real coefficients $c_{\0}$ and $c_{\1}$ by using, instead
of the second part of the Niven algorithm,  the (complex) eigenvalues of
the companion matrix associated with the quaternionic polynomial
$A_{n}(q)$. The SPV approach, avoiding the solution of the coupled real
equations (step 2), surely simplifies the Niven algorithm and gives a more
practical method to find the zeros of $n$-order unilateral quaternionic
polynomials. However, once the real coefficients $c_{\0}$ and $c_{\1}$ are
obtained in terms of the moduli and real parts of the complex eigenvalues
of the companion matrix, we have to find the quaternionic coefficients
$d_{\0}$ and $d_{\1}$. It seems that we have {\em no} alternative choice
and we have to come back to the Niven algorithm (step 1).

In this paper, we show that the matrix approach (based on the use of the
quaternion eigenvalue problem) leads us to find {\em directly} the
solutions of unilateral quaternion polynomials by calculating the
eigenvectors of the companion matrix. This allow us to {\em completely}
avoid the use of the Niven algorithm and, consequently, to hardly simplify
the method to find the zeros of unilateral quaternionic polynomials.

We interrupt our discussion at this point to introduce, the Niven
algorithm and the vector-matrix notation. These tools would permit
a more pleasant reading of the article and rend our exposition
self-consistent.

\section*{II. THE NIVEN ALGORITHM}
In this section, we aim to concentrate  on the Niven algorithm. We shall
first discuss the method proposed by Niven for a generic $n$-order
quaternionic polynomial, then, to illustrate  the potential and problems of
the Niven algorithm, we consider an explicit example, i.e. a quadratic
polynomial.

In order to determine the quaternionic coefficients  $d_{\0}$ and
$d_{\1}$ in terms of the real coefficients of the second order
polynomial $C_{\2}(q)$ and of the quaternionic coefficients of the
$n$-order polynomials $A_{n}(q)$, we expand the polynomial product
$B_{n-\2}(q)C_{\2}(q)$,
\begin{eqnarray}
\label{peb}  B_{n-\2}(q)C_{\2}(q)& =&  q^{n} - (b_{n \mi
\3}+c_{\1})\,q^{n \mi \1} - ( b_{n
\mi \4}- c_{\1}\, b_{n \mi \3} +c_{\0} )\, q^{n\mi \2} +\nonumber \\
 & &  -\sum_{s \eq \2}^{n
\mi\3}\,( b_{s
\mi \2}- c_{\1}\, b_{s \mi \1} - c_{\0}\, b_{s})\,q^s + \nonumber \\
  & & + \,(c_{\1}\,b_{\0} + c_{\0} \,b_{\1})\, q + c_{\0}\,
  b_{\0}\,\,.
\end{eqnarray}
For a particular choice of the quaternionic coefficients of the
polynomial $B_{n-\2}(q)$, i.e.
\begin{eqnarray}
\label{coeq}
b_{n \mi \3} & = & a_{n \mi \1} - c_{\1}\,\,, \nonumber \\
b_{n \mi \4}& = & a_{n \mi \2} + c_{\1}\, b_{n \mi \3} -c_{\0}\,\,, \nonumber \\
b_{s \mi \2} & = & a_{s} + c_{\1}\, b_{s \mi \1} + c_{\0}\,
b_s\,\,\,\,\,\,\,\,\,\,\mbox{\footnotesize $ (s =2,3,\cdots,n -
3)$}\,\,,
\end{eqnarray}
we can rewrite Eq.(\ref{peb}) as follows
\begin{eqnarray}
\label{e11} B_{n-\2}(q)C_{\2}(q) & = & A_{n}(q) + (a_{\1}
+c_{\1}\,b_{\0} + c_{\0} \,b_{\1})\, q + a_{\0} + c_{\0}\,
b_{\0}\,\,.
\end{eqnarray}
Comparing this  equation with the polynomial equation (\ref{pe}),
we obtain (step 1)
\begin{eqnarray}
\label{ed}
 d_{\1} & = & a_{\1}
+c_{\1}\,b_{\0}(c_{\0,\1};a_{\2,\3,\cdots,n \mi \1}) + c_{\0}
\,b_{\1}(c_{\0,\1};a_{\2,\3,\cdots,n \mi \1})\,\,,
 \nonumber \\
d_{\0}  & = &a_{\0} + c_{\0}\, b_{\0}(c_{\0,\1};a_{\2,\3,\cdots,n
\mi \1})\,\,.
\end{eqnarray}
This show that by some simple algebraic manipulations, we can express the
quaternionic coefficients $d_{\0,\1}$ in terms of the real coefficients
$c_{\0,\1}$ and of the quaternionic coefficients $a_{\0,\1,...,n\mi\1}$.
This represents  step 1 of the Niven algorithm. As stated in the
introduction, to complete the Niven algorithm, we have to determine the
real coefficients $c_{\0,\1}$ (step 2). Once such coefficients are
obtained, we can explicitly calculate the coefficients $d_{\0,\1}$  by
using Eq.(\ref{ed}) which, after step 2, only contains  known quantities.
It is natural to ask why the quaternionic coefficients $d_{\0,\1}$ are
important in calculating the quaternionic solution of the polynomial
$A_{n}(q)$. The answer is given by observing that if $q_*$ is a solution of
$A_{n}(q)$ [$\Rightarrow A_n(q_*)=0$]  and
$\{\,c_{\0}\,,\,c_{\1}\,\}=\{-\,|q_*|^{\2}\,,\,2\, \mbox{Re}[q_*]\,\}$
[$\Rightarrow C_{\2}(q_*)=0$], then  from Eq.(\ref{pe}), we find
\begin{equation}
\label{solpe} D_{\1}(q_*)=0\,\,\,\,\,\Rightarrow\,\,\,\,\, q_*= -
\,\bar{d}_{\1}\,d_{\0}\,/\,|d_{\1}|^{\2} \,\,.
\end{equation}
The quaternionic solution can thus be  expressed in terms of the
coefficients $d_{\0,\1}$. Turning to the problem of determining the real
coefficients $c_{\0,\1}$, Niven noted that
\[ c_{\0}= - \,|q_*|^{\2} = - \,|d_{\0}|^{\2}\,/\, |d_{\1}|^{\2}
\hspace*{.5cm}\mbox{and}\hspace*{.5cm} c_{\1}=\,2\, \mbox{Re}[q_*] =- 2\,
\mbox{Re}[\bar{d}_{\1}\,d_{\0}]\,/\,|d_{\1}|^{\2}\,\,.
\]
Consequently, in order to find $c_{\0,\1}$ we have to solve the
following system (step 2)
\begin{eqnarray}
\label{coupled} c_{\0} \,|d_{\1}(c_{\0,\1};a_{\1,\2,\3,\cdots,n
\mi \1})|^{\2}
+|d_{\0}(c_{\0,\1};a_{\0,\2,\3,\cdots,n \mi \1})|^{\2} & = &0\,\,,\nonumber \\
c_{\1} \,|d_{\1}(c_{\0,\1};a_{\1,\2,\3,\cdots,n \mi \1})|^{\2} +
2\,\mbox{Re} [\bar{d}_{\1}(c_{\0,\1};a_{\1,\2,\3,\cdots,n \mi
\1})d_{\0}(c_{\0,\1};a_{\0,\2,\3,\cdots,n \mi \1})] & = & 0\,\,.
\end{eqnarray}
 Each {\em real} coupled solution
$\{\,c_{\0}\,,\,c_{\1}\,\}$ gives the modulus and the real part of
the quaternionic polynomial
solution. \\

\noindent {$\bullet$ {\bf Second order quaternionic polynomials}. The main
(practical) problem in using the Niven algorithm is in step 2, that is
finding the real solutions of the coupled system (\ref{coupled}). To
illustrate the difficulty in using the Niven algorithm, we solve the
following quaternionic quadratic equation
\begin{equation}
\label{exa1}
 q^{\2} + j\,q + (1-k)=0\,\,\,\,\,[\,A_{\2}(q)=q^{\2} -a_{\1} \,q  -
 a_{\0}\,\,,\,\,\, a_{\0}=k-1\,\,,\,\,\,a_{\1}=-j\,]\,\,.
 \end{equation}
Second order polynomials represent the more simple situation in
which we can test the Niven algorithm. In this case, the
quaternionic polynomial $B_{n \mi \2}(q)$ reduces to
$B_{\0}(q)=-b_{\0}=1$. Consequently,
\begin{eqnarray*}
B_{\0}(q)C_{\2}(q)& =&  q^{\2} - c_{\1} \,q  -  c_{\0} \\
 & =  &  q^{\2} - a_{\1} \,q  -  a_{\0} + (\,a_{\1}-c_{\1}\,)\, q +
 a_{\0}-c_{\0} \\
  &=  & A_{\2}(q)+ D_{\1}(q)\,\,.
\end{eqnarray*}
From the previous equation or equivalently by setting $b_{\1}=0$
and $b_{\0}=-1$ in Eq.(\ref{ed}), we obtain
\begin{eqnarray}
\label{d2o}
d_{\1} & =&  a_{\1} -  c_{\1}\,\,,\nonumber  \\
d_{\0} & =&  a_{\0} -  c_{\0}\,\,.
\end{eqnarray}
Now, for second order polynomials, Eq.(\ref{coupled}) reduces to
\begin{eqnarray}
\label{coupled2} c_{\0} \,|a_{\1} - c_{\1}|^{\2}
+|a_{\0}- c_{\0}|^{\2} & = &0\,\,,\nonumber \\
c_{\1} \,|a_{\1} - c_{\1}|^{\2} + 2\,\mbox{Re} [(\bar{a}_{\1} -
c_{\1})(a_{\0}- c_{\0})] & = & 0\,\,.
\end{eqnarray}
The solution of this system is {\em not} simple. In this special
case, the choice of
\[a_{\1}=-j\hspace*{.5cm}\mbox{and}\hspace*{.5cm}a_{\0}=k-1\,\,,\]
reduces the system (\ref{coupled2}) to
\begin{eqnarray}
\label{coupled3} c_{\0} \,(1 + c_{\1}^{\2})+ (1 + c_{\0})^{\2} + 1
& = &0\,\,,\nonumber \\
c_{\1} \,(1 + c_{\1}^{\2})+ 2\,c_{\1}\,(1 + c_{\0}) & = &0\,\,.
\end{eqnarray}
The discussion can be then simplified by considering
\[c_{\1}\neq0\hspace*{.5cm}\mbox{and}\hspace*{.5cm}c_{\1}=
0\,\,.\] In the first case  no real coupled solution exists. For
$c_{\1}=0$, we have
\[ (\,c_{\0}\,,\,c_{\1}\,) = \left\{ \begin{array}{l} (-1,0)\\
(-2,0) \end{array}
\right.\,\,\,\,\,\,\,\,\,\,\Rightarrow\,\,\,\,\,\,\,\,\,\,
(\,d_{\0}\,,\,d_{\1}\,) = \left\{ \begin{array}{l} (k,-j)\\
(k+1,-j) \end{array} \right.\,\,\,.
\]
Finally, by using Eq.(\ref{solpe}) we conclude that the solutions of the
quadratic equation (\ref{exa1}) are given by
\begin{equation}
 q_*=\left\{ \begin{array}{l} -\,i\\
-\,(i+j) \end{array} \right.\,\,\,.
\end{equation}
By direct substitution, we can easily verify that $-i$ and $-(i+j$)
represent the zeros of our quadratic polynomial. Before concluding this
brief review of the Niven algorithm and underlining that it surely is an
important (old) technique for obtaining the solution of $n$-order
unilateral quaternionic polynomial, we wish to emphasize that the system
(\ref{coupled2}) and, in general, the system (\ref{ed}) are {\em not} very
practical and their solutions could require laborious calculations.

\section*{III. COMPANION MATRIX AND SPV METHOD}

Taking as our guide the standard complex matrix theory, we
associate to the $n$-order unilateral quaternionic polynomial
$A_n(q)$ the following companion matrix
\begin{equation}
\label{cmat}
 \left(
\begin{array}{ccccc}
a_{n \mi \1} & a_{n \mi \2} & \cdots & a_{\1} & a_{\0}\\
1 & 0 & \cdots & 0 & 0\\
0 & 1 & \cdots & 0 & 0 \\
\vdots & \vdots & \ddots & \vdots & \vdots\\
0 & 0 & \cdots & 1 & 0
\end{array}  \right)\,\,.
\end{equation}
In the standard complex theory, the zeros of $n$-order polynomials can be
determined by calculating the eigenvalues of the corresponding companion
matrix. Any direct attempt based upon a {\em straightforward} expansion of
the determinant $|M_n[\mathbb{H}]-\lambda \,I|$ is surely destined for
failure due to the non-commutativity of the quaternionic algebra
\cite{CDELA}.  The quaternionic eigenvalue problem has recently been
investigated by means of the complex representation of quaternionic
matrices.  It is clear that the rigorous presentation of the quaternionic
eigenvalue problem would require a deep knowledge   of quaternionic matrix
theory. To facilitate the understanding of this paper, we do not discuss
any theoretical aspect of the quaternionic eigenvalue problem here and aim
to give a practical matrix solution method based on complex translation
rules. We refer the reader to the articles cited in
refs.\cite{DS00,DSS02,J05} for a through discussion of quaternionic
diagonalization and Jordan form. In order to appreciate some of the
potentialities of the complex representation of quaternionic matrices, we
shall consider in detail the case of second-order polynomials and show how
we can immediately obtain the coefficients $c_{\0}$ and $c_{\1}$ by
calculating the real part and the modulus of the eigenvalues associated to
the complex matrix which represents the complex translation of the
quaternionic companion matrix (SPV method \cite{SER01}).

Before analyzing the companion matrix associated to a second-order
 unilateral quaternionic polynomial and  calculate its eigenvalues, we give
the main idea of the complex translation by using a generic  $2
\times 2$ quaternionic matrix,
\[ M_{\2}[\mathbb{H}] =
\left( \begin{array}{cc}
q_{\1} & q_{\2} \\
q_{\3} & q_{\4} \end{array} \right)\,\,,\hspace*{.5cm}
q_{\1,\2,\3,\4} \in \mathbb{H}\,\,.
\]
The complex translation  would proceed along the following lines.
A first possibility is the use of the symplectic form directly for
the quaternionic matrix $M_{\2}[\mathbb{H}]$, i.e.
\begin{equation}
\label{repc1}
 M_{\2}[\mathbb{H}] = \left( \begin{array}{cc}
z_{\1} & z_{\2} \\
z_{\3} & z_{\4} \end{array} \right) + j\,
\left( \begin{array}{cc}
w_{\1} & w_{\2} \\
w_{\3} & w_{\4} \end{array} \right)\,\,\,\rightarrow\,\,\,
M_{\4}[\mathbb{C}] = \left( \begin{array}{rrrr}
\,\,z_{\1} & \,\,\,\,z_{\2}  & - \bar{w}_{\1} &  - \bar{w}_{\2} \\
z_{\3} & z_{\4}  & - \bar{w}_{\3} &  - \bar{w}_{\4} \\
w_{\1} & w_{\2}  & \bar{z}_{\1} &  \bar{z}_{\2} \\
w_{\3} & w_{\4}  & \bar{z}_{\3} &  \bar{z}_{\4}
\end{array}
\right) \,\,,
\end{equation}
where $z_{\1,\2,\3,\4}$ and $w_{\1,\2,\3,\4} \in \mathbb{C}$. In
this case the symplectic decomposition directly applies to the
complex matrix blocks. Another possibility is represented by the
use of the symplectic decomposition for the quaternionic elements
of $M_{\2}[\mathbb{H}]$, i.e.
\begin{equation}
\label{repc2} M_{\2}[\mathbb{H}] = \left( \begin{array}{cc}
z_{\1} + j \, w_{\1} & z_{\2} + j \, w_{\2} \\
z_{\3} + j \, w_{\3} & z_{\4} + j \, w_{\4} \end{array} \right)
\,\,\,\rightarrow\,\,\, \widetilde{M}_{\4}[\mathbb{C}] = \left(
\begin{array}{rrrr}
\,\,z_{\1} & - \bar{w}_{\1} & \,\,\,\,z_{\2}   &  - \bar{w}_{\2} \\
w_{\1}  & \bar{z}_{\1} & w_{\2} &  \bar{z}_{\2} \\
z_{\3} & - \bar{w}_{\3} & z_{\4}  &   - \bar{w}_{\4} \\
w_{\3}  & \bar{z}_{\3} & w_{\4}  &  \bar{z}_{\4}
\end{array}
\right) \,\,.
\end{equation}
For the problem of finding the zeros of quaternionic polynomial,
the complex representations (\ref{repc1}) and (\ref{repc2})
represent equivalent choices. In fact,  there exist a permutation
matrix
\[P_{\4}= \left( \begin{array}{cccc}
\,1\, & \,0\, & \,0\, & \,0\, \\
0 & 0 & 1 & 0 \\
0 & 1 & 0 & 0 \\
0 & 0 & 0 & 1 \\
\end{array}
\right)
\]
which convert $M_{\4}[\mathbb{C}]$ to
$\widetilde{M}_{\4}[\mathbb{C}]$ by using a simple  similarity
transformation,
\[  \widetilde{M}_{\4}[\mathbb{C}]  = P_{\4} \, M_{\4}[\mathbb{C}]
\,P_{\4}\,\, . \] In Appendix A, we extend this similarity for the
matrices $M_{\2 n}[\mathbb{C}]$ and $\widetilde{M}_{\2
n}[\mathbb{C}]$ (complex counterparts of the $n$-dimension
quaternionic matrix $M_{n}[\mathbb{H}]$),
\[  \widetilde{M}_{\2 n}[\mathbb{C}]  = P_{\2 n} \, M_{\2
n}[\mathbb{C}]\, P_{\2 n}\,\,,\] by explicitly giving the
permutation matrix  $P_{\2 n}$. Such a matrix will be obtained
 by working with the complex representations of the
quaternionic column vectors associated to the matrix
representations (\ref{repc1}) and (\ref{repc2}). Due to this
equivalence, in all  the explicit examples given in this paper, we
shall use, without loss of generality, the complex
matrix representation $M_{\2 n}[\mathbb{C}]$. \\

\noindent {$\bullet$ {\bf Second order quaternionic polynomials}. Let us
consider the unilateral quaternionic polynomial discussed in the previous
section, Eq.(\ref{exa1}). In this case, the quaternionic companion matrix
is
\begin{equation}
\label{qm1} M_{\2}[\mathbb{H}] = \left( \begin{array}{rc}
- j   & \,\,\,k-1\,\, \\
1 & 0 \end{array} \right)\,\,.
\end{equation}
The SPV method allow to obtain the real coefficient $c_{\0}$ and
$c_{\1}$ by calculating the eigenvalues of the complex matrix
$M_{\4}[\mathbb{C}]$,
\[\{\,
\lambda_{\1}\,,\,\lambda_{\2}\,\,,\,\,\bar{\lambda}_{\1}\,,\,\bar{\lambda}_{\2}\,\}\,\,,
\] which, due to the particular structure of the complex
symplectic translation, always appear in conjugate pairs \cite{DS00}. By
using the complex translation rules (\ref{repc1}), we immediately  obtain
the following complex counterpart of $M_{\2}[\mathbb{H}]$,
\begin{equation}
\label{cm1} M_{\4}[\mathbb{C}] = \left( \begin{array}{rrrr}
0 & -1 & \,\,\,1 &-\,i\,\\
1 & 0 & 0 & 0\, \\
-1 &  - \,i &0 & -1\,\\
0 & 0 & 1 & 0\,
\end{array}
\right)\,\,,
\end{equation}
whose eigenvalues can be  easily calculated by using
$|M_{\4}[\mathbb{C}]-\lambda \,I_{\4}|=0$. As observed above, the
eigenvalues appear in conjugate pairs,
\[\{\,\lambda_{\1}=
i\,,\,\lambda_{\2}=i\,\sqrt{2}\,\,,\,\,\bar{\lambda}_{\1}=-\,i\,,\,
\bar{\lambda}_{\2}=-\,i\,\sqrt{2}\,\}\,\,.
\]
The real coefficients $c_{\0}$ and $c_{\1}$,
\[ (\,c_{\0}\,,\,c_{\1}\,) = \left\{ \begin{array}{lcl}
(\,-|\lambda_{\1}|^{\2}\,,\,2\,\mbox{Re}[\lambda_{\1}]\,) & = & (-1,0)\\
(\,-|\lambda_{\2}|^{\2}\,,\,2\,\mbox{Re}[\lambda_{\2}]\,) & = &
(-2,0) \end{array} \right.\,\,,
\]
can now be  obtained without the need  to solve the system
(\ref{coupled3}). This method proposed by Ser\^odio, Pereira and Vitoria
\cite{SER01} avoids the solution of the complicated coupled system
(\ref{coupled}) and represents a simpler way to find the real coefficients
$c_{\0}$ and $c_{\1}$. Consequently, it  improves the Niven algorithm. Note
that  step 1 in the Niven algorithm has {\em yet} to be used to find the
quaternionic coefficients $d_{\0}$ and $d_{\1}$. In the next section, we
shall prove that by calculating the eigenvectors of the complex
representation of the quaternion companion matrix $M_{n}[\mathbb{H}]$, we
can {\em directly} obtain the zeros of $n$-order unilateral quaternionic
polynomials.

\section*{IV. FINDING THE ZEROS BY USING THE EIGENVECTORS METHOD}
In the complex theory, the zeros of the $n$-order polynomial, $A_n(z)=0$,
are determined by calculating the eigenvalues of the corresponding
companion matrix, $M_n[\mathbb{C}]$, i.e. $|M_n[\mathbb{C}]-\lambda\,
I|=0$. For quaternions, this method cannot be carried through successfully.
In fact, the eigenvalues of the quaternionic companion matrix are {\em
not}\, univocally determined. Consequently, we have to pursue a different
track. Before demonstrating two theorems which will
 guide us to the solution of our quaternionic mathematical challenge,
 we briefly recall the main results concerning  the (right)
 eigenvalue problem for quaternionic matrices. Due to the
 non-commutativity of the quaternionic algebra, the eigenvalues of
\begin{equation}
\label{compmat}
 M_{n}[\mathbb{H}]  =
 \underbrace{\left(
\begin{array}{ccccc}
a_{n \mi \1}^{(\1,i)} & a_{n \mi \2}^{(\1,i)} & \cdots & a_{\1}^{(\1,i)} & a_{\0}^{(\1,i)}\\
1 & 0 & \cdots & 0 & 0\\
0 & 1 & \cdots & 0 & 0 \\
\vdots & \vdots & \ddots & \vdots & \vdots\\
0 & 0 & \cdots & 1 & 0
\end{array}  \right)}_{\displaystyle{Z_{n}[\mathbb{C}]}} +\,
j\, \underbrace{\left[-j\, \left(
\begin{array}{ccccc}
a_{n \mi \1}^{(j,k)} & a_{n \mi \2}^{(j,k)} & \cdots &
a_{\1}^{(j,k)} &
a_{\0}^{(j,k)}\\
0 & 0 & \cdots & 0 & 0\\
0 & 0 & \cdots & 0 & 0 \\
\vdots & \vdots & \ddots & \vdots & \vdots\\
0 & 0 & \cdots & 0 & 0
\end{array}  \right) \, \right]}_{\displaystyle{W_{n}[\mathbb{C}]}}
\end{equation}
belong to the following class of equivalence \begin{equation}
\{\,u_{\1}\lambda_{\1}\bar{u}_{\1}\,,\,u_{\2}\lambda_{\2}\bar{u}_{\2}\,,\,
... \, , \, u_{n}\lambda_{n}\bar{u}_{n}\,\} \,,\end{equation}
where (by using the standard definition  adopted  in literature)
$\lambda_n$
 are the eigenvalues of
\begin{equation}
\label{compmattra}
 M_{\2 n}[\mathbb{C}] =
\left(
\begin{array}{rr}
\,Z_{n}[\mathbb{C}] & \,\,\,- W^*_{n}[\mathbb{C}]\,
\\ & \\W_{n}[\mathbb{C}] & Z_{n}^*[\mathbb{C}]\,
\end{array}  \right)
\end{equation}
with $\mbox{Im}[\lambda_n] \geq 0$, and $u_n$ are quaternionic unitary
similarity transformations. The starting point in finding the zeros of
$n$-order quaternionic unilateral polynomials is the same  as in the
complex theory, i.e. the solution of the eigenvalue problem for the
companion matrix,
\begin{equation} \label{epro} M_{n}[\mathbb{H}]\,\Phi =
\Phi\,\lambda\,\,,\hspace*{.5cm}\lambda \in
\mathbb{C}\,\,,\,\,\,\mbox{Im}[\lambda] \geq 0\,\,,
\end{equation}
 and
\[ \Phi =
\left(
\begin{array}{c}
\varphi_{\1} \\ \varphi_{\2} \\ \vdots \\
\varphi_{n}
\end{array}
\right)\,\,,\hspace*{.5cm}\varphi_{\1,\2,...,n} \in
\mathbb{H}\,\,.
\]
Nevertheless, once the complex eigenvalues have been calculated we also get
 infinite {\em similar} quaternionic eigenvalues
\begin{equation} \label{epro2} M_{n}[\mathbb{H}]\,\Psi
=\Psi \,u\lambda\bar{u}\,\,,
\end{equation}
where $\Psi=\Phi\,\bar{u}$, $u\in\mathbb{H}$, and $u\bar{u}=1$.
Geometrically, the quaternionic eigenvalue $u\lambda\bar{u}$ represents,
for the imaginary part of $\lambda$, the following three-dimensional
rotation:
\[ \hspace*{.8cm} \begin{array}{c}
\lambda \\
\downarrow\\
\hspace*{-.3cm}(\mbox{Im}[\lambda],0,0) \\ \swarrow\,\searrow \\
\hspace*{-2cm} \theta=2\,\arctan
[|(\mbox{Re}[-iu],\mbox{Re}[-ju],\mbox{Re}[-ku])|/u_{\0}]\,\,\,,\,\,\,\,\,
\mbox{axis :}\,\,(\mbox{Re}[-iu],\mbox{Re}[-ju],\mbox{Re}[-ku])\\
\searrow\,\,\swarrow\\
(\mbox{Re}[-iu\lambda\bar{u}],\mbox{Re}[-ju\lambda\bar{u}],
\mbox{Re}[-ku\lambda\bar{u}])\\
\downarrow \\
\,\,\,\,u\lambda \bar{u}\,\,.\end{array}  \] We demonstrate later that, for
each complex eigenvalue, $\lambda$, it is possible to choose an appropriate
similarity transformation $u^{(q)}\lambda \bar{u}^{(q)}$ which corresponds
to one of the zeros of the $n$-order quaternionic unilateral polynomial
$A_n(q)$. This breaks down the symmetry between the infinitely many
equivalent directions in the eigenvalue quaternionic space, $u\lambda
\bar{u}$.   The complete set of the zeros of the $n$-order quaternionic
unilateral polynomial $A_n(q)$  will be then given by
\begin{equation}
\label{ss}
 \{\,u_{\1}^{(q_{\1})}\lambda_{\1}\bar{u}_{\1}^{(q_{\1})}\,,\,
 u_{\2}^{(q_{\2})}\lambda_{\2}\bar{u}_{\2}^{(q_{\2})}\,,\,
... \, , \, u_{n}^{(q_n)}\lambda_{n}\bar{u}_{n}^{(q_{n})}\,\}
\,\,.
\end{equation}
Now, we proceed with the main objective of this Section and prove that the
eigenvalue spectrum (\ref{ss}) is connected with the zeros of $A_n(q)$. In
doing that, the first step will be the proof of two theorems. Theorem A and
Theorem B concern the
(right) complex eigenvalue problem (\ref{epro}).\\

\noindent {\bf Theorem A}. {\em The last component of the
quaternionic eigenvector $\Phi$ is non null.}\\

\noindent{\bf Proof}. Due to the special form of the companion
matrix, from Eq.(\ref{epro}) we get the following equalities:
\begin{eqnarray}
\label{sis1}
 a_{n-\1}\, \varphi_{\1} + a_{n-\2}\,\varphi_{\2} +
\cdots + a_{\1}\,
\varphi_{n-\1} + a_{\0}\, \varphi_n  & = &  \varphi_{\1}\, \lambda\,\,,  \nonumber \\
\varphi_{\1}  & = &  \varphi_{\2}\, \lambda\,\,, \nonumber \\
 & \vdots &  \nonumber\\
\varphi_{n-\1} &  =  & \varphi_n \,\lambda\,\,.
\end{eqnarray}
Thus, $\varphi_n = 0$ implies $\varphi_{\1,\2,...,n-\1}=0$. This means
\[ \varphi_n = 0 \,\,\,\Rightarrow\,\,\,\Phi=\boldsymbol{0}\,\,\,\,\,
\mbox{\tiny{$\blacksquare$}}\]

\noindent {\bf Theorem B}.  {\em Let $\Phi$ be an eigenvector of
the companion matrix $M_n[\mathbb{H}]$ associated with  the
polynomial $A_n(q)$. The solution $q$, expressed in terms of the
two last components of $\Phi$,  is given by
$\varphi_{n-\1}\,\varphi^{\mi \1}_n$.}\\

\noindent{\bf Proof}. Combining the last $n-1$ equalities in
(\ref{sis1}), we obtain
\begin{eqnarray*}
\varphi_{n-\1} &  =  & \varphi_n \,\lambda\,\,,\\
\varphi_{n-\2}  & = &  \varphi_{n - \1}\, \lambda =  \varphi_n \,\lambda^{\2}\,\,,\\
 & \vdots & \\
\varphi_{\1}  & = &  \varphi_{ \2}\, \lambda =  ... = \varphi_n
\,\lambda^{n - \1} \,\,.
\end{eqnarray*}
This implies
\[\Phi = \left(
\begin{array}{c}
\varphi_{\1} \\ \varphi_{\2} \\ \vdots \\
\varphi_{n}
\end{array}
\right)=\varphi_{n}\, \left(
\begin{array}{c}
\lambda^{n \mi \1} \\  \lambda^{n \mi \1}\\ \vdots \\
1
\end{array}
\right)\,\,.
\]
The  first equation  of system (\ref{sis1}) can be then re-written
as follows
\[
 a_{n-\1}\, \varphi_{n}\,\lambda^{n-\1} + a_{n-\2}\,\varphi_{n}\,\lambda^{n-\2} +
\cdots + a_{\1} \,\varphi_{n}\,\lambda + a_{\0}\, \varphi_n   =
\varphi_{n}\, \lambda^{n}\,\,.
\]
Theorem A guarantees that $\varphi_{n}\neq 0$. We can thus
multiply (from the right) the previous equation by $\varphi^{\mi
\1}_{n}$ obtaining
\[
 a_{n-\1}\, \varphi_{n}\,\lambda^{n-\1}\,\varphi^{\mi \1}_{n} +
 a_{n-\2}\,\varphi_{n}\,\lambda^{n-\2}\,\varphi^{\mi \1}_{n}  +
\cdots + a_{\1} \,\varphi_{n}\,\lambda \,\varphi^{\mi \1}_{n} +
a_{\0}  = \varphi_{n}\, \lambda^{n}\,\varphi^{\mi \1}_{n}\,\,,
\]
or equivalently
\[
 a_{n-\1}\, (\,\varphi_{n}\,\lambda\,\varphi^{\mi \1}_{n}\,)^{n- \1} +
 a_{n-\2}\,(\,\varphi_{n}\,\lambda\,\varphi^{\mi \1}_{n}\,)^{n-\2}  +
\cdots + a_{\1} \,\varphi_{n}\,\lambda \,\varphi^{\mi \1}_{n} +
a_{\0}  = (\,\varphi_{n}\, \lambda\,\varphi^{\mi \1}_{n}\,)^n\,\,.
\]
Consequently, $q=\varphi_{n}\, \lambda\,\varphi^{\mi \1}_{n}$ and
by using that $\varphi_{n\mi \1}=\varphi_{n}\lambda$
\begin{equation}
q=\varphi_{n-\1}\, \varphi^{\mi \1}_{n}\,\,.
\end{equation}
This completes the proof.  {\tiny $\blacksquare$}\\

\noindent The results of Theorem A and Theorem B can be easily
extended to the (right) quaternionic eigenvalue problem
(\ref{epro2}). In fact,
\[ "\, \Phi\neq\boldsymbol{0} \,\,\,\,\,\Rightarrow\,\,\,\,\,\varphi_n \neq
0 \,"\,\,\,\,\, \Leftrightarrow \,\,\,\,\, "\,\Phi \, \bar{u}\,
(=\Psi) \neq\boldsymbol{0}
\,\,\,\Rightarrow\,\,\,\varphi_n\,\bar{u} \,(=\psi_n) \neq 0\,"
\]
 and
\[
q=\varphi_{n-\1}\, \varphi^{\mi \1}_{n} \,\,\,\,\,\Rightarrow\,\,\,\,\,
q=\varphi_{n-\1}\,\bar{u}\,( \varphi_{n} \,\bar{u})^{\mi
\1}=\varphi_{n-\1}\,\bar{u}\,u\, \varphi^{\mi \1}_{n}=\psi_{n-\1}\,
\psi^{\mi \1}_{n}\,\,.\] We conclude this subsection  going back to the
discussion on the quaternionic similarity transformation which breaks down
the symmetry between the infinitely many equivalent directions in the
eigenvalue quaternionic space. By using Theorem B, we can write the zero of
the quaternionic polynomial $A_n(q)$ in terms of the last component of the
eigenvector $\Phi$ and of the complex eigenvalue $\lambda$,
\begin{equation}
 q = \varphi_n \,\lambda \,\varphi_n^{\mi \1}\,\,.
\end{equation}
This allows us to define (for the companion matrix) a {\em
privileged} (right) quaternionic eigenvalue problem, i.e.
\begin{eqnarray}
\label{epro3} M_{n}[\mathbb{H}]\,\Psi^{(q)} & = &
\Psi^{(q)}\,u^{(q)}\,\lambda\,\,\bar{u}^{(q)}\hspace*{.5cm}
(u^{(q)}=\varphi_n/|\varphi_n|\,)\nonumber \\
 & = & \Psi^{(q)}\,q\,\,,
\end{eqnarray}
 with
\[ \Psi^{(q)} = \Phi\, \bar{u}^{(q)}  = \Phi\,\varphi_n^{\mi
\1}\,|\varphi_n| =  \left(
\begin{array}{c}
\varphi_{\1}\varphi_n^{\mi \1} \\ \varphi_{\2}\varphi_n^{\mi
\1} \\ \vdots \\
\varphi_{n\mi \1 }\varphi_n^{\mi \1}\\
1
\end{array}
\right)\,|\varphi_n|\,\,.
\]
Recall that
\[ q= \varphi_{n-\1}\, \varphi^{\mi \1}_{n}
\hspace*{.5cm}\mbox{and}\hspace*{.5cm} \lambda= \varphi^{\mi \1}_{n}\,
\varphi_{n-\1}\,\,.\] For the complex case, commutativity guarantees
$q=\lambda$. The problem of finding the zeros of $n$-order unilateral
quaternionic polynomials is enormously simplified by using the results
obtained in this Section. This will be explicitly seen in the next
subsections, where the solutions  of second and third order quaternionic
unilateral equations are calculated using complex matrix translation. For
the convenience of the readers, we illustrate the steps of our method in
the following sequence:
\begin{center}
\begin{tabular}{cl}
$A_n(q)$ &  ($n$-order unilateral quaternionic polynomial)\,\,,\\
 $\Downarrow $ &\\
$M_{n}[\mathbb{H}]$ & (quaternionic companion matrix)\,\,,\\
 $\Downarrow $ &\\
$M_{\2 n}[\mathbb{C}]$ & (complex translated matrix)\,\,,\\
$\Downarrow $ & \\
$(\omega_{\1},\,\omega_{\2},\, ...,\,
\omega_n,\,
\sigma_{\1},\,\sigma_{\2},\, ...,\, \sigma_n)_{\lambda} $ &
(eigenvector of $M_{\2 n}[\mathbb{C}]$  corresponding to the complex eigenvalue
 $\lambda$)\,\,,\\
 $\Downarrow $ & \\
 $(\omega_{\1}+j\,\sigma_{\1},\, \omega_{\2}+j\,\sigma_{\2},\, ...,\,
 \omega_n + j\,  \sigma_n)_{\lambda} $ &
(quaternionic translated eigenvector $\Phi_{\lambda}$)\,\,,\\
 $\Downarrow $ & \\
 $q=\varphi_{n-\1}\, \varphi^{\mi \1}_{n}$ & (zero corresponding
 to $\lambda$)\,\,,\\
 $\Downarrow $ & \\
$\lambda= \varphi^{\mi \1}_{n}\, \varphi_{n-\1}$  & (check of the
solution)\,\,.
 \end{tabular}
 \end{center}
The importance of the {\em new} method is due, in great extent, to the
possibility of using complex matrix analysis.

\subsection*{IV-A. SECOND ORDER UNILATERAL POLYNOMIALS}
In Section II and in Section III we have solved (by means of the
SPV algorithms) the following (unilateral) quaternionic quadratic
equation
\[ q^{\2} + j\, q + (1-k) =0\,\,.\]
We now calculate the zeros of $A_{\2}(q)$ by using our method:
\[
\begin{array}{rcl}
A_{\2}(q) &=& q^{\2} + j\, q + (1-k)\,\,,\\
&\Downarrow  &\\
 M_{\2}[\mathbb{H}] & = &\left( \begin{array}{rc}
- j   & \,\,\,k-1\,\, \\
1 & 0 \end{array} \right)\,\,,\\
& \Downarrow & \\
M_{\4}[\mathbb{C}] & = & \left( \begin{array}{rrrr}
0 & -1 & \,\,\,1 &-\,i\,\\
1 & 0 & 0 & 0\, \\
-1 &  - \,i &0 & -1\,\\
0 & 0 & 1 & 0\,
\end{array}
\right)\,\,,\\
& \Downarrow & \\
\left( \begin{array}{c} \,\omega_{\1}\, \\ \omega_{\2} \\
\sigma_{\1} \\ \sigma_{\2} \end{array} \right) & = & \left\{ \,
\left( \begin{array}{c} \,0\,\\ 0 \\
i \\ 1 \end{array} \right)_i\,,\,\left( \begin{array}{c} \sqrt{2}\,(\sqrt{2}-1)\\
-\,i\,(\sqrt{2}-1)\\
i\,\sqrt{2} \\ 1 \end{array} \right)_{i\,\sqrt{2}}\,\right\}\,\,,\\
& \Downarrow & \\
 \left( \begin{array}{c} \,\varphi_{\1}\, \\ \varphi_{\2}
\end{array} \right) = \left( \begin{array}{c} \,\omega_{\1} + j\, \sigma_{\1}\,
\\ \omega_{\2} + j\, \sigma_{\2}
\end{array} \right) & = &  \left\{ \,
\left( \begin{array}{r} -\,k\\ j
\end{array} \right)_{i}\,,\,\left( \begin{array}{c} \sqrt{2}\,(\sqrt{2}-1 -k) \\
j -\,i\,(\sqrt{2}-1)
\end{array} \right)_{i\,\sqrt{2}}\,\right\}\,\,,\\
& \Downarrow & \\
q= \varphi_{\1}\,  \varphi^{\mi \1}_{\2}  & = &      \{\,- \,i
\,,\, -\,(i+j)\,\}\,\,\,\left[\, \Rightarrow \,\,\,  \lambda=
\varphi^{\mi \1}_{\2}\,\varphi_{\1} =\{\,i \,,\, i\,\sqrt{2}
\,\}\,\right]\,\,.
\end{array}
\]

\subsection*{IV-B. THIRD ORDER UNILATERAL POLYNOMIALS}
Let us now give an example of solution for a third order
quaternionic  polynomial.
\[
\begin{array}{rcl}
A_{\3}(q) & = &  q^{\3} + k\, q^{\2} + i\, q - j\,\,,\\
 & \Downarrow & \\
M_{\3}[\mathbb{H}]& = & \left(\begin{array}{rrr}
-k & \,-i & \,\,\,\,j\,\,\\
1 & 0 & 0\,\,\\
0 & 1 & 0\,\, \end{array} \right)\,\,,\\
 & \Downarrow & \\
M_{\6}[\mathbb{C}] & = & \left(\begin{array}{rrrrrr}
\,0 & -i & \,\,\,0 & \,\,\,\,i & \,\,\,0 & -1\,\\
1 & 0 & 0 & 0 & 0 & 0\, \\
0  & 1& 0 & 0 &0 & 0\, \\
i  & 0&1 &0 & i& 0\,\\
0  &0 &0  & 1& 0& 0\,\\
0  &0 &0 &0 &1 & 0\,
 \end{array} \right)\,\,,\\
& \Downarrow & \\
\left( \begin{array}{c} \,\omega_{\1}\, \\ \omega_{\2} \\
\omega_{\3}\\ \sigma_{\1} \\ \sigma_{\2} \\ \sigma_{\3}\end{array}
\right) & = & \left\{\,\left(\begin{array}{r}
-1\\i\\1\\-1\\i\\1\end{array} \right)_i\,,\,\left(\begin{array}{c}
i-1\\i\sqrt{2}\\
i+1\\ i \sqrt{2}\\
i-1\\ \sqrt{2}
\end{array} \right)_{\frac{1+i}{\sqrt{2}}}\,,\,\left(\begin{array}{c}
1-i\\- \sqrt{2}\\
i+1\\ - i \sqrt{2}\\
i-1\\ \sqrt{2}
\end{array} \right)_{\frac{i-1}{\sqrt{2}}}\,
\right\}\,\,,\\
& \Downarrow & \\
\left( \begin{array}{c} \,\varphi_{\1}\, \\ \varphi_{\2} \\
\varphi_{\3}
\end{array} \right) = \left( \begin{array}{c} \,\omega_{\1} + j\, \sigma_{\1}\,
\\ \omega_{\2} + j\, \sigma_{\2} \\ \omega_{\3} + j\, \sigma_{\3}
\end{array} \right)
&=  & \left\{\,\left(\begin{array}{c} -1 -j\\i-k\\1+j\end{array}
\right)_i\,,\,\left(\begin{array}{c}
i-1-k \sqrt{2}\\
i \sqrt{2} + j-k\\
1+i+j \sqrt{2}
\end{array} \right)_{\frac{1+i}{\sqrt{2}}}\,,\,
\left(\begin{array}{c}
1-i+k\sqrt{2}\\
-\sqrt{2} - j- k\\
1+i+j\sqrt{2}
\end{array} \right)_{\frac{i-1}{\sqrt{2}}}\,
\right\}\,\,,\\
& \Downarrow & \\
 q = \varphi_{\2}\,  \varphi^{\mi \1}_{\3} & = &
\left\{\,-k\,,\, \frac{j-k+\sqrt{2}}{2}\,,\,
\frac{j-k-\sqrt{2}}{2}\,\right\}\,\,\,\left[\, \Rightarrow \,\,\,
\lambda = \varphi^{\mi \1}_{\3} \varphi_{\2}   =  \left\{\,
i\,,\,\frac{1+i}{\sqrt{2}}\,,\,\frac{i-1}{\sqrt{2}}\,\right\}\,\right]\,\,.
\end{array}
\]

\section*{V. BILATERAL QUADRATIC EQUATIONS}
So far, we have not analyzed the possibility to have left and
right acting coefficients in the quaternionic polynomials. This
topic exceeds the scope of this paper. At the moment it represents
an open question and  satisfactory resolution methods  have to be
further investigated. Nevertheless, an interesting generalization
of the matrix approach introduced in the previous section can be
found in discussing bilateral equations\cite{FAPESP}, i.e.
\begin{equation}
\label{ebi} p^{\2} - \alpha_{\1}\,p + p\,\beta_{\1} - \alpha_{\0}
= 0\,\,.
\end{equation}
For second order bilateral quaternionic polynomial,  we introduced
the following {\em generalized}  companion matrix
\begin{equation}
{\mathcal M}_{\2}[\mathbb{H}] = \left(\begin{array}{cc} \,\alpha_{\1}\, & \,
\alpha_{\0}\,\\
1 & \beta_{\1} \end{array} \right)\,\,,
\end{equation}
which reduces to the standard one for $\beta_{\1}=0$. The (right)
complex eigenvalue problem for this matrix,
\begin{equation} \label{eprob} {\mathcal M}_{\2}[\mathbb{H}]\,\Phi =
\Phi\,\lambda \,\,,
\end{equation}
is equivalent to
\begin{eqnarray*}
\alpha_{\1}\,\varphi_{\1} + \alpha_{\0}\,\varphi_{\2} & = &
\varphi_{\1}\,\lambda\,\,,\\
\varphi_{\1} + \beta_{\1}\,\varphi_{\2} & = &
\varphi_{\2}\,\lambda\,\,,
\end{eqnarray*}
where $\varphi_{\1,\2}$ are the components of the eigenvector
$\Phi$. It is immediate to show that $\varphi_{\2}=0$ implies
$\Phi=0$. Consequently, the second component  of the (non trivial)
eigenvector $\Phi$ is always different from zero. Multiplying the
first of the previous equations  by $\varphi_{\2}^{\mi \1}$ (from
the right), we obtain
\[\alpha_{\1}\,\varphi_{\1}\,\varphi_{\2}^{\mi
\1} + \alpha_{\0}  = \varphi_{\1}\,\varphi_{\2}^{\mi
\1}\,\varphi_{\2}\,\lambda\,\varphi_{\2}^{\mi \1}\,\,.
\]
From the second equation, we get
\[
\varphi_{\2}\,\lambda\,\varphi_{\2}^{\mi
\1}=\varphi_{\1}\,\varphi_{\2}^{\mi \1} + \beta_{\1}\,\,.
\]
Finally,
\[
\alpha_{\1}\,\varphi_{\1}\,\varphi_{\2}^{\mi \1} + \alpha_{\0}  =
(\varphi_{\1}\,\varphi_{\2}^{\mi \1})^{^{\2}} +
\varphi_{\1}\,\varphi_{\2}^{\mi \1} \, \beta_{\1}
\]
which, once compared  with Eq.(\ref{ebi}), implies
\begin{equation}
p=\varphi_{\1}\,\varphi_{\2}^{\mi\1}\,\,.
\end{equation}
The solution for the bilateral quadratic equation (\ref{ebi}) is
then {\em formally} equal to the solution given for the unilateral
one. The complex eigenvalue $\lambda$ can be expressed in terms of
the $\Phi$ components and of the right acting coefficient
$\beta_{\1}$ by
\begin{equation}
\lambda=\varphi_{\2}^{\mi \1}\,\varphi_{\1} + \varphi_{\2}^{\mi
\1}\, \beta_{\1}\, \varphi_{\2}\,\,.
\end{equation}
It is worth pointing that the bilateral equation (\ref{ebi}) can
be reduced to the following {\em equivalent} unilateral equation
\begin{equation}
q^{\2} - a_{\1}\,q - a_{\0} = 0\,\,,
\end{equation}
where $a_{\1} =  \alpha_{\1}
 + \beta_{\1}$ and
$a_{\0}  =  \alpha_{\0} - \alpha_{\1}\,\beta_{\1}$ by using
\begin{equation}
p=q-\beta_{\1}\,\,. \end{equation} The next examples should emphasize one
more time the advantage of using  of the matrix approach proposed in this
paper to solve unilateral (and some particular bilateral) quaternionic
equations. It would be desirable to extend this approach to more general
polynomials (with left and right acting quaternionic coefficients) but we
have not been able to do this and from our point of view this should
deserve more attention and further investigation. The question of solving
general quaternionic polynomials is at present far from being solved.

\subsection*{V-A. SECOND ORDER BILATERAL POLYNOMIALS}
Let us consider the bilateral quadratic equation
\begin{equation}
\label{bil} p^{\2} - i\, p + p \, j -k =0\,\,.
\end{equation}
To solve this equation, we shall follow the same steps of the
resolution method presented in the previous section:
\[
\begin{array}{rcl}
p^{\2} - i\, p &+& p \, j -k\,\,,\\
&\Downarrow  &\\
{\mathcal M}_{\2}[\mathbb{H}] & = &\left( \begin{array}{cc}
\,\,i\,\,   & \,\,k\,\, \\
1 & j \end{array} \right)\,\,,\\
& \Downarrow & \\
{\mathcal M}_{\4}[\mathbb{C}] & = & \left( \begin{array}{rrrr}
\,i & 0 & \,\,\,0 &\,-\,i\,\\
\,1 & 0 & 0 & -\,1\, \\
0 &  - \,i &\,-\,i & 0\,\\
0 & 1 & 1 & 0\,
\end{array}
\right)\,\,,\\
& \Downarrow & \\
\left( \begin{array}{c} \,\omega_{\1}\, \\ \omega_{\2} \\
\sigma_{\1} \\ \sigma_{\2} \end{array} \right) & = & \left\{ \,
\left( \begin{array}{c} \,1\,\\ 0 \\
0 \\ 1 \end{array} \right)_{\0}\,,\,\left( \begin{array}{c} -1-\sqrt{2}\\
i+i\,\sqrt{2}\\
-\,i \\ 1 \end{array} \right)_{i\,\sqrt{2}}\,\right\}\,\,,\\
& \Downarrow & \\
 \left( \begin{array}{c} \,\varphi_{\1}\, \\ \varphi_{\2}
\end{array} \right) = \left( \begin{array}{c} \,\omega_{\1} + j\, \sigma_{\1}\,
\\ \omega_{\2} + j\, \sigma_{\2}
\end{array} \right) & = &  \left\{ \,
\left( \begin{array}{c} 1\\ j
\end{array} \right)_{\0}\,,\,\left( \begin{array}{c} k-1-\sqrt{2}\\
i + i\,\sqrt{2}+j
\end{array} \right)_{i\,\sqrt{2}}\,\right\}\,\,,\\
& \Downarrow & \\
p= \varphi_{\1}\,  \varphi^{\mi \1}_{\2}  & = &      \{\,- \,j
\,,\, i\,\}\,\,\,\left[\, \Rightarrow \,\,\,  \lambda=
\varphi^{\mi \1}_{\2}\,\varphi_{\1} +  \varphi^{\mi
\1}_{\2}\,\beta_{\1}\, \varphi_{\2}  =\{\,0 \,,\, i\,\sqrt{2}
\,\}\,\right]\,\,.
\end{array}
\]

\subsection*{V-B. EQUIVALENT SECOND ORDER UNILATERAL POLYNOMIALS}
In this subsection, we solve the unilateral quadratic  equation
obtained from Eq.(\ref{bil}) by using $p=q-j$:
\[
\begin{array}{rcl}
q^{\2} &-& (i+j)\, q \,\,,\\
&\Downarrow  &\\
 M_{\2}[\mathbb{H}] & = &\left( \begin{array}{cc}
i+ j   & \,\,\,0\,\, \\
1 & 0 \end{array} \right)\,\,,\\
& \Downarrow & \\
M_{\4}[\mathbb{C}] & = & \left( \begin{array}{rrrr}
\,i & \,\,\,0 & \,-\,1 &\,\,\,0\,\\
\,1 & 0 & 0 & 0\, \\
1 &  0 &-\,i & 0\,\\
0 & 0 & 1 & 0\,
\end{array}
\right)\,\,,\\
& \Downarrow & \\
\left( \begin{array}{c} \,\omega_{\1}\, \\ \omega_{\2} \\
\sigma_{\1} \\ \sigma_{\2} \end{array} \right) & = & \left\{ \,
\left( \begin{array}{c} \,0\,\\ 1 \\
0 \\ 1 \end{array} \right)_{\0}\,,\,\left( \begin{array}{c} -\,\sqrt{2}\,(1+\sqrt{2})\\
i\,(1+\sqrt{2})\\
i\,\sqrt{2} \\ 1 \end{array} \right)_{i\,\sqrt{2}}\,\right\}\,\,,\\
& \Downarrow & \\
 \left( \begin{array}{c} \,\varphi_{\1}\, \\ \varphi_{\2}
\end{array} \right) = \left( \begin{array}{c} \,\omega_{\1} + j\, \sigma_{\1}\,
\\ \omega_{\2} + j\, \sigma_{\2}
\end{array} \right) & = &  \left\{ \,
\left( \begin{array}{c} 0\\ 1+j
\end{array} \right)_{i}\,,\,\left( \begin{array}{c} - \sqrt{2}\,(1+ \sqrt{2}+k) \\
i\,(1+ \sqrt{2})+j
\end{array} \right)_{i\,\sqrt{2}}\,\right\}\,\,,\\
& \Downarrow & \\
q= \varphi_{\1}\,  \varphi^{\mi \1}_{\2}  & = &      \{\,0 \,,\,
i+j\,\}\,\,\,\left[\, \Rightarrow \,\,\,  \lambda= \varphi^{\mi
\1}_{\2}\,\varphi_{\1} =\{\,0 \,,\, i\,\sqrt{2} \,\}\,\right]\,\,.
\end{array}
\]
As expected, by recalling that
\[ p = q -\beta_{\1} = q - j =  \{\,-\,j \,,\,
i\,\}\] we recover the solution of Eq.(\ref{bil}) obtained in the
previous subsection.

\section*{VI. CONCLUSIONS}

In this paper we have presented a practical method to find the roots of
unilateral quaternion polynomials by means of the eigenvalue spectra and
eigenvectors of the companion matrix. In particular, we have shown that the
zeros of unilateral quaternion polynomials can be directly given in terms
of the two last component of the eigenvectors of the companion matrix. The
practical resolution method is enormously simplified by using the complex
matrix translation. We have also extended the method to bilateral equations
introducing the {\em generalized} companion matrix. One question still
unanswered is whether similar approaches could be used to finding the zeros
of $n$-order polynomial with left and right acting quaternionic
coefficients. The paper is essentially self-contained and the possibility
to translate the quaternionic eigenvalue problem to the corresponding
complex one, then, come back to quaternions, and give the final solutions
should represent a fundamental tool for the reader which is not so familiar
with the quaternionic algebra. In this spirit, the applications (explicit
examples of $2$ and $3$ order quaternionic polynomials) would clarify  the
advantage in adopting this {\em full} matrix approach  in solving
quaternionic unilateral  equations.\\

\noindent {\bf ACKNOWLEDGEMENTS.} The authors wish to express their
gratitude to Prof. Nir Cohen for several helpful comments and for many
stimulating conversations. One of the author (S.D.L) greatly thanks the
Department of Mathematics (Parana University), where the paper was
prepared, for the invitation and hospitality. This is part of the research
project financed by the Araucaria Foundation (V.L and G.D).

\newpage

\section*{APPENDIX A. EQUIVALENCE BETWEEN TRANSLATION RULES}
Let
\begin{equation}
\label{vc1} \left(
\begin{array}{c}
\varphi_{\1} \\ \varphi_{\2} \\ \vdots \\
\varphi_{n}
\end{array}
\right)  =  \left(
\begin{array}{c}
\omega_{\1}\\  \omega_{\2}\\ \vdots \\
\omega_{n}
\end{array}
\right) + j\, \left(
\begin{array}{c}
\sigma_{\1}\\  \sigma_{\2}\\ \vdots \\
\sigma_{n}
\end{array}
\right) \to \left(
\begin{array}{c}
\omega_{\1}\\  \omega_{\2}\\ \vdots \\
\omega_{n}\\ \sigma_{\1}\\  \sigma_{\2}\\ \vdots \\
\sigma_{n}
\end{array}
\right)
\end{equation}
and
\begin{equation}
\label{vc2} \left(
\begin{array}{c}
\varphi_{\1} \\ \varphi_{\2} \\ \vdots \\
\varphi_{n}
\end{array}
\right)  =  \left(
\begin{array}{c}
\omega_{\1} + j \, \sigma_{\1} \\  \omega_{\2} + j\,\sigma_{\2} \\ \vdots \\
\omega_{n} + j \, \sigma_n
\end{array}
\right)
 \to \left(
\begin{array}{c}
\omega_{\1}\\ \sigma_{\1}\\ \omega_{\2}\\ \sigma_{\2}\\ \vdots \\\vdots \\
\omega_{n}\\ \sigma_{n}
\end{array}
\right)
\end{equation}
be the complex vectors obtained by using  the complex translation
rules given in  Eqs.(\ref{repc1}) and (\ref{repc2}). A simple
algebraic calculation shows that
\begin{equation}
\left(
\begin{array}{c}
\omega_{\1}\\ \sigma_{\1}\\ \omega_{\2}\\ \sigma_{\2}\\ \vdots \\\vdots \\
\omega_{n}\\ \sigma_{n}
\end{array}
\right) \, = \,P_{\2 n}\,\left(
\begin{array}{c}
\omega_{\1}\\  \omega_{\2}\\ \vdots \\
\omega_{n}\\ \sigma_{\1}\\  \sigma_{\2}\\ \vdots \\
\sigma_{n}
\end{array}
\right)\,\,,
\end{equation}
where
\begin{equation}
P_{\2 n\,,\,rs} = \left\{ \begin{array}{cl} 1 & \hspace*{.5cm}
(r,s)=(1,1)\,,\,\,(2,n+1)\,,\,\,(3,2)\,,\,\,(4,n+2)\,,\,\, \dots
\,,\,\,
(2n-1,n)\,,\,\,(2n,2n)\,\,,\\
0 &  \hspace*{.5cm} \mbox{otherwise}\,\,.
\end{array}
\right.
\end{equation}
The foregoing result can be then used to prove the equivalence
between the complex matrices $M_{\2 n}[\mathbb{C}]$ and
$\widetilde{M}_{\2 n}[\mathbb{C}]$,
\begin{equation}
\widetilde{M}_{\2 n}[\mathbb{C}]  = P_{\2 n} \, M_{\2
n}[\mathbb{C}] \,P_{\2 n}\,\, .
\end{equation}

\end{document}